\documentclass[preprint,11pt]{elsarticle}

\usepackage{amssymb}
\usepackage{amsmath}
\usepackage{amsfonts}
\usepackage{amssymb}
\usepackage{indentfirst,latexsym,bm}
\usepackage{amsthm}
\usepackage{color,xcolor}
\usepackage[all]{xy}
\input{mathrsfs.sty}
\newtheorem{theorem}{Theorem}[section]
\newtheorem{lemma}[theorem]{Lemma}
\newtheorem{corollary}[theorem]{Corollary}
\newtheorem{proposition}[theorem]{Proposition}

\newtheorem{definition}[theorem]{Definition}
\newtheorem{remark}[theorem]{Remar}

\theoremstyle{question}

\theoremstyle{problem}

\numberwithin{equation}{section}

\journal{XXX}

\begin{document}

\begin{frontmatter}



% \linenumbers


% \title{A supplement to: ``Numerical radius in Hilbert $C^*$-modules" [Math. Inequal. Appl. 24 (2021), no. 4, 1017--1030]}

\title{Notes on the numerical radius for adjointable operators on Hilbert $C^*$-modules}
\author[shnu]{Jia Li}
\ead{1184048619@qq.com}
\author[shnu]{Kangjian Wu}
\ead{wukjcool@163.com}
\author[shnu]{Qingxiang Xu}
\ead{qingxiang\_xu@126.com}
%\cortext[cor2]{Corresponding author}
%\fntext[fn1]{Partially supported by the
%National Natural Science Foundation of China (11971136).}
\address[shnu]{Department of Mathematics, Shanghai Normal University, Shanghai 200234, PR China}
\begin{abstract}Given a Hilbert module $H$ over a $C^*$-algebra, let $\mathcal{L}(H)$ be the set of all adjointable operators on $H$.
  For each $T\in\mathcal{L}(H)$, its numerical radius  is defined by $w(T)=\sup\big\{\|\langle Tx, x \rangle\|: x\in H, \|x\|=1\big\}$.
  It is proved that $w(T)=\|T\|$ whenever $T$ is normal. Examples are constructed to show that there exist Hilbert module $H$ over certain $C^*$-algebra and $T_1,T_2\in \mathcal{L}(H)$ with
  $T_1^2=0$ such that $w(T_1)\ne \frac12 \|T_1\|$ and
  $\sup\limits_{\theta\in [0,2\pi]}\|\mbox{Re}(e^{i\theta}T_2)\|<w(T_2)$.
  In addition, a new characterization of the spatial numerical radius is given, and it is proved  that
  $w\big(\pi(T)\big)\le w(T)$ for every faithful representation  $(\pi, X)$ of $\mathcal{L}(H)$ and every $T\in\mathcal{L}(H)$.
  Some inequalities are derived based on the newly obtained results.
    % Based on some elementary propositions, it is shown that some results in the paper indicated by the title can be directly derived as well as be refined.

  % Based on some elementary propositions, it is shown that some results about numerical radius in the \cite{Zamani: Hilbert module numerical radius} and papers \cite{CL} can be  directly derived as well as be refined.

\end{abstract}


\end{frontmatter}

% \linenumbers


\section{Introduction}\label{sec:intro}

Let $\mathfrak{A}$ be a general $C^*$-algebra. Recall that a linear space $H$ is said to be a (right) Hilbert $\mathfrak{A}$-module if $H$ is a (right) $\mathfrak{A}$-module and there is an $\mathfrak{A}$-valued inner-product $\langle x,y\rangle$ for each pair $x$ and $y$ in $H$ such that $H$ is complete with respect to the induced norm $$\|x\|=\sqrt{\|\langle x,x\rangle\|} \quad (x\in H).$$
For the standard references on Hilbert $C^*$-modules, the reader is referred to \cite{Lance,MT,Paschke}. Note that the $C^*$-algebra $\mathfrak{A}$ can be regarded as a Hilbert module over itself in the  usual way:
$$\langle a,b\rangle=a^*b\quad (a,b\in \mathfrak{A}).$$
Note also that every Hilbert module over the complex field $\mathbb{C}$ is exactly a Hilbert space whose inner-product is linear with respect to the second variable. So, Hilbert $C^*$-modules are natural generalizations of $C^*$-algebras as well as Hilbert spaces.

There are many dissimilarities between Hilbert spaces and Hilbert $C^*$-modules. For instance, when $x$ and $y$ are two elements in a Hilbert module over a general $C^*$-algebra such that $x\perp y$ (i.e.,  $\langle x,y\rangle=0$), it may happen that
$\|x+y\|^2\ne \|x\|^2+\|y\|^2$. A closed submodule $M$ of a Hilbert $C^*$-module $H$ may fail to be orthogonally complemented, that is, it may happen that $H\ne M+M^\perp$. Also, many new phenomena may happen when one deals with operators on Hilbert $C^*$-modules compared with that on Hilbert spaces.

Suppose that $\mathfrak{A}$ is a $C^*$-algebra, $H$ is a Hilbert $\mathfrak{A}$-module. Denote by  $\mathcal{L}(H)$ the set of all adjointable operators on $H$. Let $\mathcal{L}(H)_{\mbox{sa}}$  be the set of all self-adjoint elements in the $C^*$-algebra $\mathcal{L}(H)$.  For each $A\in\mathcal{L}(H)$, let $\mathcal{R}(A)$ and $\mathcal{N}(A)$ be  its range and  null space, respectively. In the special case that $H$ is a Hilbert space, the notation $\mathbb{B}(H)$ is used instead of $\mathcal{L}(H)$. When $b$ is an element of a unital $C^*$-algebra $\mathfrak{B}$, let $\sigma_{\mathfrak{B}}(b)$ denote the spectrum of $b$ with respect to $\mathfrak{B}$.

Unless otherwise specified, throughout the rest of this paper, $\mathfrak{A}$ is a $C^*$-algebra, $\mathcal{S}(\mathfrak{B})$ is the set consisting of all states on $\mathfrak{B}$, and $H$ is a  Hilbert $\mathfrak{A}$-module.
For each $T\in \mathcal{L}(H)$, its numerical radius is defined by
\begin{equation}\label{equ:defn of nr} w(T)=\sup\big\{\|\langle Tx, x \rangle\|: x\in H, \|x\|=1\big\}.
\end{equation}
It is known (see e.g.\,\cite[P.\,78 and Theorem~2.2]{MH} and \cite[Lemma~1]{HM}) that
\begin{align}&\label{equ:between half and one}\frac12 \|T\|\le w(T)\le \|T\|,\quad \forall T\in\mathcal{L}(H),\\
\label{equ:self-adjoint}&w(T)=\|T\|\quad\mbox{whenever}\quad T\in\mathcal{L}(H)_{\mbox{sa}}.
\end{align}
In the special case that $H$ is a Hilbert space, the following  results are well-known (see e.g.\,\cite[Theorems~2.2.1 and 2.2.11]{BDMP}):
\begin{align}\label{equality-01}&w(T)=\|T\|\quad \mbox{if $T\in \mathbb{B}(H)$ is normal},\\
\label{equality-02}&w(T)=\sup_{\theta\in [0,2\pi]}\|\mbox{Re}(e^{i\theta}T)\|,\quad \forall T\in\mathbb{B}(H),\\
\label{equality-03}&w(T)=\frac12 \|T\|\quad\mbox{if $T\in\mathbb{B}(H)$ satisfies $T^2=0$.}
\end{align}

A huge number of papers focusing on the numerical radius can be found in the literature. However,  as far as we know little has been done  concerning the validity of
\eqref{equality-01}--\eqref{equality-03} in the  Hilbert $C^*$-module case. The main purpose of this paper is to give the positive answer to
the validity of \eqref{equality-01}, meanwhile to construct examples to show that both of \eqref{equality-02} and \eqref{equality-03}
may fail to be true in the Hilbert $C^*$-module case; see Corollary~\ref{cor:normal is ok} and Theorem~\ref{thm:counterexample} for the details.

The rest of this paper is organized as follows. The spatial numerical radius is dealt with in Section~\ref{sec:snr}.
The validity of \eqref{equality-01} and examples of invalidity of \eqref{equality-02}--\eqref{equality-03} are investigated in Section~\ref{sec:validity} for Hilbert $C^*$-module operators. In addition, a new characterization of the spatial numerical radius is given in
Theorem~\ref{thm:key-01}, and it is proved in Corollary~\ref{cor:rep is less} that
$w\big(\pi(T)\big)\le w(T)$ for every faithful representation  $(\pi, X)$ of $\mathcal{L}(H)$ and every $T\in\mathcal{L}(H)$.
Based on Theorem~\ref{thm:key-01} and Corollary~\ref{cor:rep is less}, some inequalities are derived in Section~\ref{sec:applications}.

\section{The spatial numerical radius}\label{sec:snr}

To get the positive answer to the validity of \eqref{equality-01}, we need
to make use of $\widetilde{w}(T)$ introduced in \cite{MM} for each $T\in\mathcal{L}(H)$ .

\begin{definition}{\rm {\cite[Definition~2.4]{MM}} For each $T\in\mathcal{L}(H)$, let
\begin{equation}\label{equ:defn of snr}\widetilde{w}(T)=\sup_{\substack{x\in H, \rho\in\mathcal{S}(\mathfrak{A})\\ \rho\langle x,x\rangle=1}}|\rho\langle x,Tx\rangle|,
\end{equation}
which is called the spatial numerical radius of $T$ because of Theorem~\ref{thm:key-01}.
}\end{definition}

We begin with a new proof of a known result.
\begin{lemma}[{cf.\,\cite[Theorem~2.13(a) and Corollary~2.14]{MM}}]\label{lem:self case}
   For every $T\in\mathcal{L}(H)$, we have
$\widetilde{w}(T)\le \|T\|$. Moreover, $\widetilde{w}(T)=\|T\|$ whenever $T$ is self-adjoint.
\end{lemma}
\begin{proof}Given any $x\in H$, it is clear that
\begin{equation*}\langle Tx,Tx\rangle=\langle T^*Tx,x\rangle\le \langle \|T^*T\|x,x\rangle=\|T\|^2\langle x,x\rangle,
\end{equation*}
which yields
\begin{equation*}\rho \langle Tx,Tx\rangle \le \|T\|^2\rho\langle x,x\rangle,\quad \forall \rho\in\mathcal{S}(\mathfrak{A}).
\end{equation*}
So for each $x\in H$ and $\rho\in \mathcal{S}(\mathfrak{A})$, we have
\begin{align*}|\rho\langle x,Tx\rangle|^2\le \rho\langle x,x\rangle\cdot \rho\langle Tx,Tx\rangle\le \|T\|^2 \left(\rho\langle x,x\rangle\right)^2.
\end{align*}
Hence
\begin{equation}\label{equ:for less that norm T}|\rho\langle x,Tx\rangle|\le \|T\|\cdot \rho\langle x,x\rangle,
\end{equation}
which leads by \eqref{equ:defn of snr} to $\widetilde{w}(T)\le \|T\|$.

Suppose now that $T=T^*$. In this case, we prove that  $\widetilde{w}(T)=\|T\|$. In view of \eqref{equ:defn of nr}--\eqref{equ:self-adjoint}, it needs only to verify that
\begin{equation}\label{equ:reduced-01}\|\langle x,Tx\rangle\|\le \widetilde{w}(T),\quad \forall x\in H, \|x\|=1.
\end{equation}
If $\langle x,Tx\rangle=0$, then the above inequality is trivially satisfied. Assume that $x\in H$ with $\|x\|=1$ such that $\langle x,Tx\rangle\ne 0$.
Since $\langle x,Tx\rangle$ is a self-adjoint element of $\mathfrak{A}$, it is known (see e.g.\,\cite[Theorem~3.3.6]{Murphy}) that there exists $\rho_0\in\mathcal{S}(\mathfrak{A})$ such that $|\rho_0\langle x,Tx\rangle|=\|\langle x,Tx\rangle\|$. This together with  \eqref{equ:for less that norm T}  yields
$$1=\|x\|^2\ge \rho_0 \langle x,x\rangle\ge \frac{\|\langle x,Tx\rangle\|}{\|T\|}>0.$$
Let $\lambda=\sqrt{\rho_0\langle x,x\rangle}$ and $u=\frac{1}{\lambda}x$. Then
$$0<\lambda\le 1,\quad \rho_0\langle u,u\rangle=1.$$
Hence,
\begin{align*}\widetilde{w}(T)\ge |\rho_0\langle u,Tu\rangle|=\frac{|\rho_0\langle x,Tx\rangle|}{\lambda^2}=\frac{\|\langle x,Tx\rangle\|}{\lambda^2}\ge \|\langle x,Tx\rangle\|.
\end{align*}
This completes the proof of \eqref{equ:reduced-01}.
\end{proof}

Next, we make a useful observation. Let $I$ and $J$ be two non-empty sets, and let
$E=\{x_{ij}:i\in I, j\in J\}$ be a  non-empty subset of the real line $\mathbb{R}$. Suppose that $E$ is bounded above. Then it is easy to verify that
$$\sup_{i\in I}\sup_{j\in J}x_{ij}=\sup_{j\in J}\sup_{i\in I}x_{ij}.$$

Now, we provide the technical result of this section  as follows.

\begin{lemma}\label{lem:key equ}{\rm \cite[Lemma~3.4]{MM}} For every $T\in\mathcal{L}(H)$, we have
\begin{equation}\label{equ:key equ}\widetilde{w}(T)=\sup_{\theta\in [0,2\pi]}\|\mbox{Re}(e^{i\theta}T)\|.
\end{equation}
\end{lemma}
\begin{proof}Let $\rho$ be arbitrary in $\mathcal{S}(\mathfrak{A})$. Note that $\rho(b)\in \mathbb{R}$ for every $b\in \mathfrak{A}_{\mbox{sa}}$,  so $$\mbox{Re}\big(\rho(a)\big)=\rho\big(\mbox{Re}(a)\big),\quad \forall a\in \mathfrak{A}.$$
It follows that for each $x\in H$ with $\rho\langle x,x\rangle=1$,
\begin{align*}|\rho\langle x,Tx\rangle|=&\sup_{\theta\in [0,2\pi]} \left|\mbox{Re}\big[e^{i\theta}\rho\langle x,Tx\rangle\big]\right|=\sup_{\theta\in [0,2\pi]} \left|\mbox{Re}\big[\rho\big\langle x,e^{i\theta} Tx\big\rangle\big]\right|\\
=&\sup_{\theta\in [0,2\pi]} \left|\rho\big\langle x,\mbox{Re}(e^{i\theta} T)x\big\rangle\right|.\end{align*}
Consequently, by Lemma~\ref{lem:self case} we have
\begin{align*}\widetilde{w}(T)=&\sup_{\substack{x\in H, \rho\in\mathcal{S}(\mathfrak{A})\\ \rho\langle x,x\rangle=1}}|\rho\langle x,Tx\rangle|=\sup_{\substack{x\in H, \rho\in\mathcal{S}(\mathfrak{A})\\ \rho\langle x,x\rangle=1}}\sup_{\theta\in [0,2\pi]} \left|\rho\langle x,\mbox{Re}(e^{i\theta} T)x\rangle\right|\\
=&\sup_{\theta\in [0,2\pi]}\sup_{\substack{x\in H, \rho\in\mathcal{S}(\mathfrak{A})\\ \rho\langle x,x\rangle=1}}\left|\rho\langle x,\mbox{Re}(e^{i\theta} T)x\rangle\right|
=\sup_{\theta\in [0,2\pi]} \widetilde{w}\big[\mbox{Re}(e^{i\theta} T)\big]\\=&\sup_{\theta\in [0,2\pi]}\|\mbox{Re}(e^{i\theta} T)\|.
\end{align*}
This shows the validity of \eqref{equ:key equ}.
\end{proof}

\section{Investigations of \eqref{equality-01}--\eqref{equality-03} in the Hilbert $C^*$-module case}\label{sec:validity}
\subsection{The validity of \eqref{equality-01} in the Hilbert $C^*$-module case}
We start with the following comparison theorem.

\begin{theorem}\label{thm:comparison}
For every $T\in\mathcal{L}(H)$, we have $\widetilde{w}(T)\le w(T)$.
\end{theorem}

\begin{proof} Given any $S\in\mathcal{L}(H)$ and $x\in H$, let
$$a=\langle \mbox{Re}(S)x,x\rangle,\quad b=\langle \mbox{Im}(S)x,x\rangle.$$
Then $a,b\in \mathfrak{A}_{\mbox{sa}}$, so there exists $\rho\in \mathcal{S}(\mathfrak{A})$ such that
$|\rho(a)|=\|a\|$.  Since $\rho(a)$ and $\rho(b)$ are real numbers, we have
\begin{align*}\|\langle Sx,x\rangle\|=\|a-ib\|\ge |\rho(a-ib)|=\sqrt{\rho(a)^2+\rho(b)^2}\ge |\rho(a)|=\|a\|.
\end{align*}
Therefore, $$\|\langle \mbox{Re}(S)x,x\rangle\|\le \|\langle Sx,x\rangle\|,\quad\forall S\in\mathcal{L}(H), x\in H.$$
It follows that for every $T\in\mathcal{L}(H)$ and $\theta\in [0,2\pi]$,
\begin{align*}
  w\big[\mbox{Re}(e^{i\theta}T)\big] & =\sup_{\|x\|=1}\big\|\langle \mbox{Re}(e^{i\theta}T)x, x\rangle\big\|\leq \sup_{\|x\|=1}\big\|\langle e^{i\theta}Tx, x\rangle\big\| \\
  &=\sup_{\|x\|=1}\|\langle Tx, x\rangle\|=w(T).
\end{align*}
So a combination of \eqref{equ:key equ}  and \eqref{equ:self-adjoint} yields
$$\widetilde{w}(T)=\sup_{\theta\in [0,2\pi]}\|\mbox{Re}(e^{i\theta}T)\|=\sup_{\theta\in [0,2\pi]} w\big[\mbox{Re}(e^{i\theta}T)\big]\leq w(T).\qedhere$$
\end{proof}

Recall that  a pair $(\pi, X)$ is said to be a representation of $\mathcal{L}(H)$, if $X$ is a Hilbert space and $\pi:\mathcal{L}(H)\to\mathbb{B}(X)$ is a $C^*$-morphism.

\begin{theorem}\label{thm:key-01} For every faithful representation  $(\pi, X)$ of $\mathcal{L}(H)$, we have
 \begin{equation*}\widetilde{w}(T)=w\big(\pi(T)\big),\quad\forall T\in\mathcal{L}(H).
 \end{equation*}
 \end{theorem}
\begin{proof} Given any $S\in\mathcal{L}(H)$, it is easily seen that
$\pi\big(\mbox{Re}(S)\big)=\mbox{Re}\big(\pi(S)\big)$. Hence, for every $T\in\mathcal{L}(H)$ and $\theta\in [0,2\pi]$ we have
\begin{align*}\big\|\mbox{Re}(e^{i\theta}T)\big\|=&\big\|\pi\big(\mbox{Re}(e^{i\theta}T)\big)\big\|=\big\|\mbox{Re}\big(\pi(e^{i\theta}T)\big)\big\|=\big\|\mbox{Re}\big(e^{i\theta}\pi(T)\big)\big\|.
\end{align*}
So, we may combine \eqref{equ:key equ}  and \eqref{equality-02} to conclude that
\begin{align*}\widetilde{w}(T)=&\sup_{\theta\in [0,2\pi]}\|\mbox{Re}(e^{i\theta}T)\|=\sup_{\theta\in [0,2\pi]}\big\|\mbox{Re}\big(e^{i\theta}\pi(T)\big)\big\|=w\big(\pi(T)\big). \qedhere
\end{align*}
\end{proof}

As applications of Theorems~\ref{thm:comparison} and \ref{thm:key-01}, we have the following useful corollaries.

\begin{corollary}\label{cor:rep is less}For every faithful representation  $(\pi, X)$ of $\mathcal{L}(H)$ and every $T\in\mathcal{L}(H)$, we have
$w\big(\pi(T)\big)\le w(T)$.
\end{corollary}
\begin{proof}It is immediate from Theorems~\ref{thm:comparison} and \ref{thm:key-01}.
\end{proof}

\begin{corollary}{\rm {\cite[Theorems~2.13 and 2.17]{MM}}\label{cor:sample}} For every $T\in\mathcal{L}(H)$, we have
\begin{equation*}\label{equ:basic inequs-copy}\frac12 \|T\|\le \widetilde{w}(T)\le \|T\|.
\end{equation*}
Furthermore, $\widetilde{w}(T)=\|T\|$ if $T$ is normal, while $\widetilde{w}(T)=\frac12 \|T\|$ if $T^2=0$.
\end{corollary}
\begin{proof}Choose any faithful representation $(\pi, X)$ of $\mathcal{L}(H)$. Given any $T\in\mathcal{L}(H)$, by Theorem~\ref{thm:key-01} and \eqref{equ:between half and one} we have
$$\widetilde{w}(T)=w\big(\pi(T)\big)\le \|\pi(T)\|=\|T\|.$$
See also Lemma~\ref{lem:self case} for the inequality $\widetilde{w}(T)\le \|T\|$. Similarly,
$$\widetilde{w}(T)=w\big(\pi(T)\big)\ge \frac12 \|\pi(T)\|=\frac12 \|T\|.$$
If $T$ is normal, then so is $\pi(T)$. Hence by \eqref{equality-01},
$$\widetilde{w}(T)= w\big(\pi(T)\big)=\|\pi(T)\|=\|T\|.$$
Similar reasoning shows that $\widetilde{w}(T)=\frac12 \|T\|$ whenever $T^2=0$.
\end{proof}

\begin{remark}{\rm To get the equality $\widetilde{w}(T)=\frac12 \|T\|$, it is assumed in \cite[Theorem~2.17]{MM} that $T^2=0$ and $\mathcal{R}(T)$ is closed in $H$. The proof of the preceding corollary indicates that the assumption of the closedness of $\mathcal{R}(T)$  is reductant.
}\end{remark}

\begin{corollary}\label{cor:normal is ok}If $T\in\mathcal{L}(H)$ is normal, then $w(T)=\|T\|$.
\end{corollary}
\begin{proof}Suppose that $T\in\mathcal{L}(H)$ is normal. By Theorem~\ref{thm:comparison} and Corollary~\ref{cor:sample}, we have
$$w(T)\ge \widetilde{w}(T)=\|T\|,$$
which is combined with \eqref{equ:between half and one} to conclude that $w(T)=\|T\|$.
\end{proof}

\subsection{Examples of invalidity of \eqref{equality-02}--\eqref{equality-03} in the Hilbert $C^*$-module case}
Our aim here is  to construct examples to show the invalidity of \eqref{equality-02} and \eqref{equality-03} in the Hilbert $C^*$-module case. For this, we need a couple of lemmas.

\begin{lemma}\label{lem:unital Cstar-alg}
Let $\mathfrak{A}$ be  a unital $C^*$-algebra regarded as the Hilbert module over itself in the usual way. Then for every $T\in \mathcal{L}(\mathfrak{A})$,
we have $w(T)=\|T\|$.
\end{lemma}
\begin{proof}Given an arbitrary $T\in \mathcal{L}(\mathfrak{A})$, by \eqref{equ:between half and one} we need only to prove that $w(T)\ge \|T\|$. Since $\mathfrak{A}$ has a unital (denoted by $e$),
there exists a unique element $a\in \mathfrak{A}$ such that $T=L_a$, where $L_a$ is defined by
\begin{equation}\label{defn of L a}L_a(b)=ab\quad (b\in\mathfrak{A}).
\end{equation}
It follows that
\begin{align*}
  w(T) =&w(L_a)=\sup\left\{\big\| \langle L_a(b), b\rangle \big\|: b\in \mathfrak{A}, \|b\|=1\right\}\\
  \ge& \big\| \langle L_a(e), e\rangle \big\|=\|a^*\|=\|a\|=\|L_a\|=\|T\|.
\end{align*}
This completes the proof.
\end{proof}

\begin{lemma}\label{lem:less}
Let $\mathfrak{A}$ be  a unital $C^*$-algebra regarded as the Hilbert module over itself in the usual way. Then for every $a\in \mathfrak{A}$, we have
\begin{equation}\label{equ:1st inequ}\sup_{\theta\in [0,2\pi]}\big\|\mbox{Re}(e^{i\theta}L_a)\big\|\le \sup\big\{|\rho(a)|:\rho\in\mathcal{S}(\mathfrak{A})\big\},
\end{equation}
in which $L_a$ is defined by \eqref{defn of L a}.
\end{lemma}
\begin{proof}Given $a\in \mathfrak{A}$ and $\theta\in [0,2\pi]$, let $T=L_a$ and $b=e^{i\theta}a$. Then
$$ \mbox{Re}(e^{i\theta}T)=\frac{1}{2}\left[e^{i\theta}T+(e^{i\theta}T)^*\right]=\frac{1}{2}(L_b+L_{b^*})
=L_{\mbox{Re}(b)}.$$
Since  $\mbox{Re}(b)$ is self-adjoint, by \cite[Theorem~3.3.6]{Murphy} we have
\begin{align*}\|\mbox{Re}(b)\|=&\sup\big\{\big|\rho[\mbox{Re}(b)]\big|:\rho\in\mathcal{S}(\mathfrak{A})\big\}=\sup\big\{\big|\mbox{Re}[\rho(b)]\big|:\rho\in\mathcal{S}(\mathfrak{A})\big\}\\
\le&\sup\big\{|\rho(b)|:\rho\in\mathcal{S}(\mathfrak{A})\big\}=\sup\big\{|\rho(a)|:\rho\in\mathcal{S}(\mathfrak{A})\big\}.
\end{align*}
Consequently,
\begin{equation*}\big\|\mbox{Re}(e^{i\theta}T)\big\|=\big\|L_{\mbox{Re}(b)}\big\|=\|\mbox{Re}(b)\|\le \sup\big\{\big|\rho(a)|:\rho\in\mathcal{S}(\mathfrak{A})\big\}.
\end{equation*}
The desired conclusion follows from the arbitrariness of $\theta$ in $[0,2\pi]$.
\end{proof}

\begin{lemma}\label{lem:special}Let  $M_2(\mathbb{C})$ be the $C^*$-algebra consisting of all $2\times 2$ complex matrices. Put
$$\mathfrak{A}=M_2(\mathbb{C}), \quad a=\left(
                                          \begin{array}{cc}
                                            1 & 0 \\
                                            1 & 0 \\
                                          \end{array}
                                        \right).$$
Then
\begin{equation}\label{equ:norm a bigger}\sup \big\{|\rho(a)|: \rho\in S(\mathfrak{A})\big\}<\|a\|.
\end{equation}
\end{lemma}
\begin{proof}Let $I_2$ be the identity matrix in $\mathfrak{A}$, and
\begin{equation}\label{equ:defn of b}e_1=\left(
        \begin{array}{cc}
          1 & 0 \\
          0 & 0 \\
        \end{array}\right),\quad e_2=I_2-e_1,\quad b=\left(
                                          \begin{array}{cc}
                                            0 & 0 \\
                                            1 & 0 \\
                                          \end{array}
                                        \right).
\end{equation}
Then for every $\rho\in \mathcal{S}(\mathfrak{A})$, we have $\rho(e_1)+\rho(e_2)=\rho(I_2)=\|\rho\|=1$. In view of $a=e_1+b$, we obtain
\begin{align*}|\rho(a)|\le \rho(e_1)+|\rho(b)|\le \rho(e_1)+\sqrt{\rho(bb^*)}=\rho(e_1)+\sqrt{\rho(e_2)}.
\end{align*}
Since $0\le \rho(e_1)\le 1$, it follows that
\begin{align*}|\rho(a)|\le \max_{0\le x\le 1} \left(x+\sqrt{1-x}\,\right)=\frac54<\sqrt{2}=\|a\|.
\end{align*}
So, the desired inequality follows.
\end{proof}

We are now in the position to provide the main result of this subsection as follows.

\begin{theorem}\label{thm:counterexample}There exist Hilbert module $H$ over certain $C^*$-algebra $\mathfrak{A}$ and $T_1,T_2\in \mathcal{L}(H)$ with
$T_1^2=0$ such that $w(T_1)\ne \frac12 \|T_1\|$ and
\begin{equation*}\sup_{\theta\in [0,2\pi]}\|\mbox{Re}(e^{i\theta}T_2)\|<w(T_2).
\end{equation*}
\end{theorem}
\begin{proof}
 Let $\mathfrak{A}$ and $a\in \mathfrak{A}$ be specified as in Lemma~\ref{lem:special} such that \eqref{equ:norm a bigger} is satisfied.
Let $H=\mathfrak{A}$, which is regarded as Hilbert module over $\mathfrak{A}$ in the usual way, and put $T_1=L_b$ and $T_2=L_a$, where
$b$ is given by \eqref{equ:defn of b}. Then $T_1^2=0$ and $T_1\ne 0$, so by Lemma~\ref{lem:unital Cstar-alg} $w(T_1)=\|T_1\|\ne \frac12 \|T_1\|$. Utilizing
\eqref{equ:1st inequ} and  \eqref{equ:norm a bigger} together with Lemma~\ref{lem:unital Cstar-alg}  yields
                             $$\sup_{\theta\in [0,2\pi]}\|\mbox{Re}(e^{i\theta}T_2)\|\le \sup\big\{\big|\rho(a)|:\rho\in\mathcal{S}(\mathfrak{A})\big\}<\|a\|=\|T_2\|=w(T_2).$$
This completes the proof.
\end{proof}

\section{Some applications}\label{sec:applications}

It is remarkable that many results can be derived by a simple use of  Theorem~\ref{thm:key-01} and Corollary~\ref{cor:rep is less}. For this, we take two propositions as examples. Our first proposition is concerned with a new proof of a known result.

\begin{proposition}\label{prop:sample}{\rm {\cite[Theorem~3.2]{MM}}} For every $T\in\mathcal{L}(H)$, we have
$$\frac14\|T^*T+TT^*\|\le \widetilde{w}^2(T)\le \frac12 \|T^*T+TT^*\|.$$
\end{proposition}
\begin{proof}Choose any faithful representation $(\pi, X)$ of $\mathcal{L}(H)$. By \cite[Theorem~1]{Kittaneh} we have
$$\frac14\|\pi(T)^*\pi(T)+\pi(T)\pi(T)^*\|\le w^2\big(\pi(T)\big)\le \frac12 \|\pi(T)^*\pi(T)+\pi(T)\pi(T)^*\|,$$
which clearly leads by Theorem~\ref{thm:key-01} to the desired inequality.
\end{proof}

\begin{remark}{\rm
 Following the line in the proof of Proposition~\ref{prop:sample}, it is easily seen that most results in \cite[Section~3]{MM} for Hilbert $C^*$-module operators can be derived directly from Theorem~\ref{thm:key-01}.
}\end{remark}

Our next proposition illustrates the usefulness of Corollary~\ref{cor:rep is less}.
\begin{proposition}For every $A_i,B_i\in\mathcal{L}(H)$ $(i=1,2)$, we have
\begin{equation*}\label{equ:mixed inequality}r(A_1B_1+A_2B_2)\le \frac12\big(w(B_1A_1)+w(B_2A_2)+\sqrt{\alpha}\,\big),
\end{equation*}
where $r(A_1B_1+A_2B_2)$ denotes the spectral radius of $A_1B_1+A_2B_2$, and
\begin{equation*}\alpha=\big[w(B_1A_1)-w(B_2A_2)\big]^2+4 \|B_1A_2\|\cdot \|B_2A_1\|.
\end{equation*}
\end{proposition}
\begin{proof} For simplicity, we put
\begin{align*}&T=A_1B_1+A_2B_2,\quad a=\|B_1A_2\|\cdot \|B_2A_1\|,\\
&f(x,y)=\frac12\left[x+y+\sqrt{(x-y)^2+4a}\right] \quad (0\le x,y<+\infty).
\end{align*}
Choose any faithful representation $(\pi, X)$ of $\mathcal{L}(H)$. Since
$$\sigma_{\mathcal{L}(H)}(T)\cup \{0\}=\sigma_{\pi[\mathcal{L}(H)]}\big(\pi(T)\big)\cup\{0\}=\sigma_{\mathbb{B}(H)}\big(\pi(T)\big)\cup\{0\},$$
we have $r(T)=r\big(\pi(T)\big)$. In virtue of $\|\pi(B_1A_2)\|\cdot \|\pi(B_2A_1)\|=a$, it can be concluded by \cite[Theorem~2.2]{AF}
that
\begin{align*}r\big(\pi(T)\big)\le& f\Big(w\big(\pi(B_1A_1)\big), w\big(\pi(B_2A_2)\big)\Big).
\end{align*}
Our aim here is to prove that
$$r(T)\le f\Big(w(B_1A_1), w(B_2A_2)\Big).$$
Due to Corollary~\ref{cor:rep is less} and $r(T)=r\big(\pi(T)\big)$, it needs only to show that
\begin{equation}\label{equ:increasment}f(x_1,y_1)\le f(x_2,y_2)\quad \mbox{whenever}\quad 0\le x_1\le x_2,0\le y_1\le y_2.
\end{equation}
If $a=0$, then $f(x,y)=\max\{x,y\}$, hence \eqref{equ:increasment} is trivially satisfied.
Suppose now that $a>0$. Then for every $x,y\in [0,+\infty)$,
\begin{align*}&\frac{\partial f}{\partial x}(x,y)=\frac{\sqrt{(x-y)^2+4a}+(x-y)}{2\sqrt{(x-y)^2+4a}}>0,\\
&\frac{\partial f}{\partial y}(x,y)=\frac{\sqrt{(x-y)^2+4a}+(y-x)}{2\sqrt{(x-y)^2+4a}}>0.
\end{align*}
So,  in this case \eqref{equ:increasment} is also valid. This completes the proof.
\end{proof}

\vspace{2ex}

%\noindent\textbf{Data Availability Statement} No dataset was generated or analyzed during this study.

%\vspace{2ex}

%\noindent\textbf{Declarations}

%\noindent\textbf{Authors' contributions} Authors declare that they have contributed equally to this paper. All authors have read and approved this version.

%\noindent\textbf{Conflict of interest} The authors have no competing interests to declare that are relevant to the content of this article.

%\begin{acknowledgements}\label{ackref}
%\end{acknowledgements}

% Important: Do not put any empty line here.

%
\end{document}